\documentclass[titlepage,11pt]{article}
\usepackage{latexsym}
\usepackage{amsfonts}
\usepackage{amsmath}
\usepackage{url}
\usepackage[T1]{fontenc}
\usepackage{enumerate}
\newcommand\blackslug{\hbox{\hskip 1pt \vrule width 4pt height 8pt depth 1.5pt
        \hskip 1pt}}
\newcommand\bbox{\hfill \quad \blackslug \bigbreak}

\def\ll{,\ldots,}

\usepackage[absolute]{textpos}
\usepackage{graphicx}

\title{Caterpillars in Erd\H{o}s-Hajnal}
\author{
Anita Liebenau\thanks{Previously at Monash University. Research supported by a DECRA Fellowship from the Australian Research Council. 
The author would like to thank for its hospitality the Institute of
Informatics, University of Warsaw, where this work was carried out.}
\\
School of Mathematics and Statistics, UNSW Sydney, NSW 2052, Australia.\\Email: a.liebenau@unsw.edu.au
\\
\\
Marcin Pilipczuk\thanks{The research of Marcin
Pilipczuk is a part of projects that have received funding from the European Research Council (ERC) under the European
Union's Horizon 2020 research and innovation programme under grant agreement No 714704.}
\\Institute of Informatics, University of Warsaw, Poland
\\
\\
Paul Seymour\thanks{Supported by ONR grant N00014-14-1-0084 and NSF
grant DMS-1265563.} and Sophie Spirkl \\
Princeton University, Princeton, NJ 08544, USA}

\date{October 18, 2017; revised \today}

\newtheorem{thm}{}[section]

\newcommand{\Proof}{\noindent{\bf Proof.}\ \ }

\begin{document}
\begin{textblock}{20}(11.3, 13.75)
\includegraphics[width=50px]{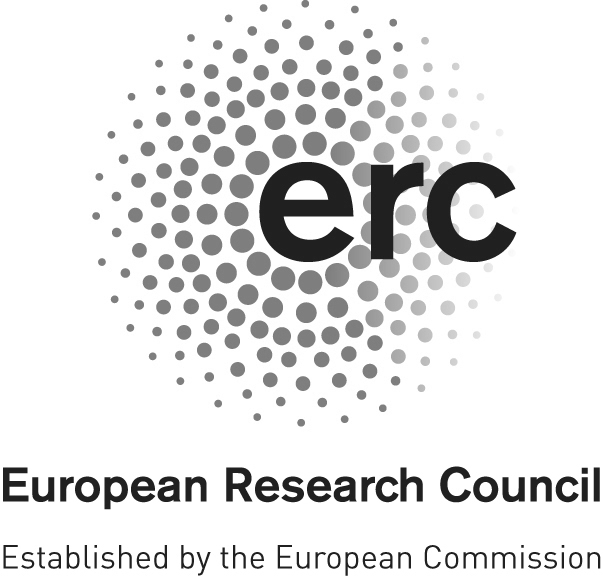}%
\end{textblock}
\begin{textblock}{20}(12.8, 13.9)
\includegraphics[width=40px]{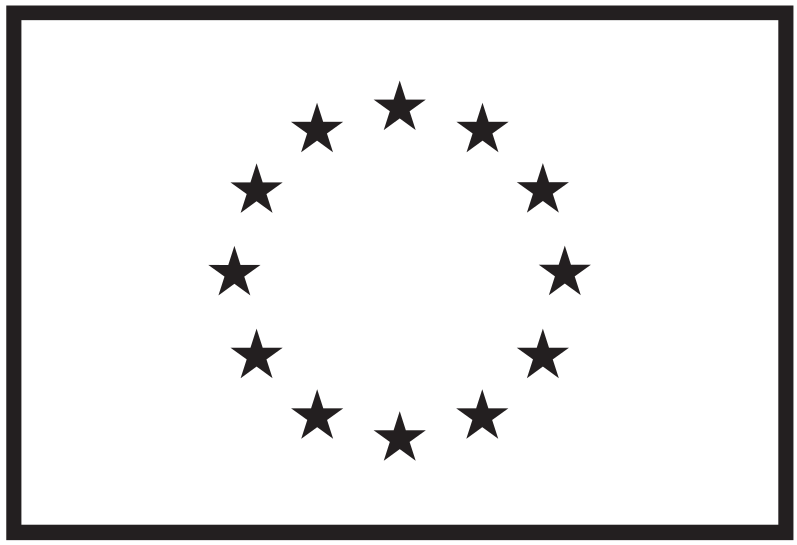}%
\end{textblock}

\maketitle

\begin{abstract}
Let $T$ be a tree such that all its vertices of degree more than two lie on one path; that is, $T$ is a 
caterpillar subdivision. We prove that there exists $\epsilon>0$ such that for every graph $G$ with $|V(G)|\ge 2$ not containing $T$ as an induced subgraph,
either some vertex has at least $\epsilon|V(G)|$ neighbours, or there are two disjoint sets of vertices $A,B$, 
both of cardinality
at least $\epsilon|V(G)|$,  where there is no edge joining $A$ and $B$. 

A consequence is: for every 
caterpillar subdivision $T$, there exists $c>0$ such that for every graph $G$ containing 
neither of $T$ and its complement as an 
induced subgraph, $G$ has a clique or stable set with at least $|V(G)|^c$ vertices. This extends a theorem of
Bousquet, Lagoutte and Thomass\'e~\cite{lagoutte}, who proved the same when $T$ is a path, and a recent theorem of 
Choromanski, Falik, Liebenau, Patel and Pilipczuk~\cite{hooks}, who proved it when $T$ is a ``hook''.
\end{abstract}

\section{Introduction}

The Erd\H{o}s-Hajnal conjecture~\cite{EH0, EH} asserts:
\begin{thm}\label{EHconj}
{\bf Conjecture: }For every graph $H$, there exists $c>0$ such that every $H$-free graph $G$ satisfies 
$$\max(\omega(G),\alpha(G))\ge |V(G)|^c.$$
\end{thm}
(All graphs in this paper are finite and have no loops or parallel edges. A graph $G$ is {\em $H$-free} 
if no induced subgraph of 
$G$ is isomorphic to $H$; and $\omega(G), \alpha(G)$ denote the cardinalities of the largest cliques and stable sets in $G$
respectively, and $\omega(G)$ is called the {\em clique number} of $G$.)
This conjecture has been investigated heavily, and nevertheless has been proved only for 
very restricted graphs $H$ (see~\cite{maria}
for a survey, and see~\cite{fox} for progress on the conjecture in a geometric setting). 
In particular it has not yet been proved when $H$ is a five-vertex path.

On the other hand, a theorem of Bousquet, Lagoutte and Thomass\'e~\cite{lagoutte} asserts the following ($\overline{H}$
denotes the complement of a graph $H$):
\begin{thm}\label{path+antipath}
For every path $H$, there exists $c>0$ such that every graph $G$ that is both $H$-free and $\overline{H}$-free
satisfies 
$\max(\omega(G),\alpha(G))\ge |V(G)|^c$.
\end{thm}
Let us say $H$ is
a {\em hook} if $H$ is a tree obtained from a path by adding a vertex adjacent to the third vertex 
of the path.
Two of us, with Choromanski, Falik, and Patel~\cite{hooks}, extended \ref{path+antipath}, proving:

\begin{thm}\label{hook+antihook}
For every hook $H$, there exists $c>0$ such that every graph $G$ that is both $H$-free and $\overline{H}$-free
satisfies
$\max(\omega(G),\alpha(G))\ge |V(G)|^c$.
\end{thm}

The main step of the proof of \ref{path+antipath} is the following:
\begin{thm}\label{getpath}
For every path $H$, there exists $\epsilon>0$ such that for every $H$-free graph $G$ with $|V(G)|\ge 2$, either 
some vertex has at least $\epsilon|V(G)|$ neighbours, or there are two anticomplete sets of vertices $A,B$,
both of cardinality
at least $\epsilon|V(G)|$.
\end{thm}
(Two sets $A,B\subseteq V(G)$ are {\em complete} to each other if $A\cap B=\emptyset$ and every vertex in $A$
is adjacent to every vertex in $B$; and {\em anticomplete} to each other if they are complete to each other in $\overline{G}$.)

It is natural to ask, which other graphs $H$ have the property of \ref{getpath}? Let us say
a graph $H$ has the {\em sparse strong EH-property} if there exists $\epsilon>0$ such that for 
every $H$-free graph $G$ with $|V(G)|\ge 2$, either
some vertex has at least $\epsilon|V(G)|$ neighbours, or there are two anticomplete sets of vertices $A,B$,
both of cardinality
at least $\epsilon|V(G)|$. Which graphs have the sparse strong EH-property?

And here is a related question: let us say a graph has the {\em symmetric strong EH-property} if there exists $\epsilon>0$ such that for
every graph $G$ that is both $H$-free and $\overline{H}$-free, with $|V(G)|\ge 2$, 
there are two disjoint sets of vertices,
both of cardinality
at least $\epsilon|V(G)|$, and either complete or anticomplete to each other. Which graphs have the symmetric strong EH-property?

It follows from a theorem of R\"odl~\cite{rodl} (and see~\cite{fox0} for a version with much better constants)
that every graph with the sparse property has the symmetric property; and 
Erd\H{o}s's construction~\cite{girth} of a graph with large girth and large chromatic number also
shows that every graph with the sparse property is a forest, 
and every graph with the symmetric property is either a forest or the complement of one. (We omit all these proofs, which
are easy; see~\cite{hooks} for more details.)
We conjecture the converses, that is:
\begin{thm}\label{treeconj}
{\bf Conjectures: } 
\begin{itemize}
\item A graph $H$ has the sparse strong EH-property if and only if $H$ is a forest.
\item A graph $H$ has the symmetric strong EH-property if and only if one of $H$, $\overline{H}$ is a forest.
\end{itemize}
\end{thm}
The first implies the second, because of the theorem of R\"odl~\cite{rodl}.
These two conjectures are reminiscent of the Gy\'arf\'as-Sumner conjecture, which we discuss later.
(Since this paper was submitted for publication, both of these conjectures have been proved, in~\cite{trees}.)

A graph $H$ is a {\em caterpillar} if $H$ is a tree and some path of $H$ contains all vertices with degree at least two;
and a {\em caterpillar subdivision} if $H$ is a tree
and some path of $H$ contains all vertices with degree at least three. 
(Thus a graph is a caterpillar subdivision if and only if it 
can be obtained from a caterpillar by subdividing edges.) We will prove:

\begin{thm}\label{mainthm}
Every caterpillar subdivision has the sparse strong EH-property.
\end{thm}
\ref{mainthm} implies the next result, which generalizes \ref{path+antipath} and \ref{hook+antihook}.
(This theorem was proved independently by the first two authors and by the last two, but since the proofs
were virtually identical we have combined the two papers into one. The original paper by the first two authors
is available~\cite{liebenau}.)
If $X\subseteq V(G)$, $G[X]$ denotes the subgraph of $G$ induced on $X$.

\begin{thm}\label{cat+anticat}
Let $H, J$ be caterpillar subdivisions. Then there exists $c>0$ such that for every graph $G$, if $G$ is both
$H$-free and $\overline{J}$-free, then 
$\max(\omega(G),\alpha(G))\ge |V(G)|^c$.
\end{thm}
\noindent{\bf Proof of \ref{cat+anticat}, assuming \ref{mainthm}.} 
There is a caterpillar subdivision such that both
$H,J$ are induced subgraphs of it, and so, by replacing $H,J$ by this graph,
we may assume that $H=J$.
Let $\epsilon$ satisfy \ref{mainthm}; so $0\le \epsilon\le 1$. By a theorem of R\"odl~\cite{rodl}, 
\\
\\
(1) {\em There exists $\delta > 0$ such that for every $H$-free graph
$G$, there is a subset $X\subseteq V(G)$ with $|X|\ge \delta |V(G)|$ such that
one of $G[X], \overline{G}[X]$ has at most $\epsilon|X|^2/4$ edges.}
\\

Choose $c$ such that $2(\epsilon\delta/2)^{2c}= 1$. 
A graph is {\em perfect} if chromatic number equals clique number for all its induced subgraphs. 
For a graph $G$, let $\pi(G)$ denote the maximum
cardinality of a subset $X$ such that $G[X]$ is perfect; we will prove by induction on $|V(G|$
that if $G$ is both
$H$-free and $\overline{H}$-free, then $\pi(G)\ge |V(G)|^{2c}$ (and consequently the theorem will follow,
since $\alpha(G)\omega(G)\ge \pi(G)$).
If $|V(G)|\le 1$ the result is trivial, and if $2\le |V(G)|\le 2/\delta$ then
$\pi(G)\ge 2\ge |V(G)|^{2c}$ as required, since $(2/\delta)^{2c} \le (2/(\epsilon\delta))^{2c} =2$.
Thus we may assume that $|V(G)|>2/\delta$.
By (1) there is a subset $X\subseteq V(G)$ with $|X|\ge \delta |V(G)|$ such that
one of $G[X], \overline{G}[X]$ has at most $\epsilon|X|^2/4$ edges; and by replacing $G$ by its complement 
if necessary, we may assume that $G[X]$ has at most $\epsilon|X|^2/4$ edges. Choose distinct 
$v_1\ll v_k\in X$, maximal such that for $1\le i\le k$, $v_i$ has at least $\epsilon|X|/2$ neighbours in
$X\setminus \{v_1\ll v_i\}$. Let $Y=X\setminus \{v_1\ll v_k\}$.
It follows that $k\le |X|/2$, and every vertex in $Y$ has fewer than $\epsilon|X|/2$ neighbours in $Y$,
from the maximality of $k$. Thus $|Y|\ge |X|/2\ge \delta|V(G)|/2$, and $G[Y]$ has maximum degree less than
$\epsilon|X|/2\le \epsilon|Y|$. 
Since $|V(G)|>2/\delta$,
it follows that $|X|>2$ and so $|Y|>1$.
By \ref{mainthm} applied to $G[Y]$, 
there are two anticomplete sets of vertices $A,B$,
both of cardinality
at least $\epsilon|Y|$. From the inductive hypothesis $\pi(G[A])\ge |A|^{2c}$ and $\pi(G[B])\ge |B|^{2c}$, and so
$$\pi(G)\ge |A|^{2c}+|B|^{2c}\ge 2 (\epsilon|Y|)^{2c}\ge 2 (\epsilon\delta|V(G)|/2)^{2c} = |V(G)|^{2c}.$$
This proves \ref{cat+anticat}.~\bbox

Let $G$ be a graph and for every subset $X\subseteq V(G)$ let $\mu(X)$ be a real number, satisfying:
\begin{itemize}
\item $\mu(\emptyset)=0$ and $\mu(V(G))=1$, and $\mu(X)\le \mu(Y)$ for all $X,Y$ with $X\subseteq Y$; and
\item $\mu(X\cup Y)\le \mu(X)+\mu(Y)$ for all disjoint sets $X,Y$.
\end{itemize}
We call such a function $\mu$ a {\em mass} on $G$. For instance, we could take $\mu(X)=|X|/|V(G)|$, or
$\mu(X) = \chi(G[X])/\chi(G)$, where $\chi$ denotes chromatic number. We denote by $N(v)$ the set of neighbours
of $v$. The result \ref{mainthm} can be extended to graphs
with masses, in the following way:

\begin{thm}\label{measureversion}
For every caterpillar subdivision $H$, there exists $\epsilon>0$ such that for every $H$-free graph $G$, and 
mass $\mu$ on $G$, either
\begin{itemize}
\item $\mu(\{v\})\ge \epsilon$ for some vertex $v$; or
\item $\mu(N(v))\ge \epsilon$ for some vertex $v$; or
\item there are two anticomplete sets of vertices $A,B$, where $\mu(A), \mu(B)\ge \epsilon$.
\end{itemize}
\end{thm}
We prove this in the next section.
It implies \ref{mainthm}, setting $\mu(X)=|X|/|V(G)|$. To see this, observe that if only the first outcome holds,
and $\mu(\{v\})\ge \epsilon$ for some $v$, then $v$ has no neighbours (or else the second outcome would hold),
and $\mu(V(G)\setminus \{v\})<\epsilon$ (or else the third outcome would hold); and so $\mu(v)>1-\epsilon$. 
Adding, $2\mu(v)> \epsilon + (1-\epsilon) = 1$, and so $\mu(v)>1/2$, and hence $|V(G)|=1$.

But \ref{measureversion} also has an interesting application to the Gy\'arf\'as-Sumner conjecture~\cite{gyarfastree,sumner}, 
which states that
for every tree $T$ and every integer $k\ge 0$, there exists $f(T,k)$ such that every $T$-free graph with clique
number at most $k$ has chromatic number at most $f(T,k)$. This has not been proved in general, and not even for
caterpillars; and not even for trees with exactly two vertices of degree more than two 
(such a tree is a simple kind of caterpillar subdivision).
But by induction on $k$, one could assume that for every vertex $v$, the chromatic number of the subgraph induced
on $N(v)$ is bounded; and so the following consequence of
\ref{measureversion} might be of interest.

\begin{thm}\label{GYconj}
Let $T$ be a caterpillar subdivision, and $k\ge 0$ an integer. Let $\epsilon$ satisfy \ref{measureversion}.
Suppose that every $T$-free graph with clique number $<k$
has chromatic number at most $c\ge 1$. Then in every $T$-free graph with clique number at most $k$ and chromatic number
more than $c/\epsilon$,  there are two anticomplete sets of vertices $A,B$, where $\chi(G[A]), \chi(G[B])\ge \epsilon \chi(G)$.
\end{thm}
\Proof
Let $G$ be a $T$-free graph with $\omega(G)\le k$.
Define $\mu(X)=\chi(G[X])/\chi(G)$, for each $X\subseteq V(G)$.
Thus one of the three outcomes of \ref{measureversion} holds.
The first implies that $\chi(G)\le 1/\epsilon$, and the second implies that $\chi(G) \le c/\epsilon$, in both
cases a contradiction.
So the third holds. This proves \ref{GYconj}.~\bbox

Incidentally, perhaps one can unify the Gy\'arf\'as-Sumner conjecture and \ref{treeconj}, in the natural way (using masses).

\section{The main proof}

In this section we prove \ref{measureversion}, but before the details of the proof, let us sketch the idea.
If $X,Y$ are disjoint subsets of $V(G)$, we say that $X$ {\em covers} $Y$ if every vertex in $Y$ has a neighbour in $X$.
First let $T$ be a caterpillar, rather than
a caterpillar subdivision, and suppose that $G$ is a $T$-free graph with a mass that does not satisfy the theorem.
We choose some large number (depending on $T$)
of disjoint subsets of $V(G)$, each with large mass (let us call them ``blocks''). 
It follows from the falsity of the third bullet of \ref{measureversion} that for every two blocks,
most of the vertices in one will have neighbours in the other, so we are well-equipped with edges between blocks. 
Choose a block $B_1$, and let us
grow a subset $X$ of it, one vertex at a time, until there is some other block, say $B_2$, 
that is at least half covered by $X$. We cannot use $B_1$ as a block any more, and we discard it, retaining only the
set $X$.
Also we discard from $B_2$ the part of $B_2$ that is not covered by $X$, and for every other block $B_3$ say, discard from $B_3$
the part that is covered by $X$.
We now have many disjoint blocks (one fewer than before), 
all still with large mass (about half what it was before), together with one more set 
$X$ that covers one of our blocks and has no edges to the others. Now pick another block (which could be $B_2$)
and do it again, 
growing a subset of it until it covers half of a different block, and so on. We can construct more complicated patterns 
of covering, by judiciously choosing which block to grow within next. 
This will enable us to find a copy of the caterpillar $T$, with all its vertices in different
blocks. 

In the case when $T$ is a caterpillar subdivision, we were not able to prove that there is 
a copy of $T$ with all its vertices in 
different blocks. But $T$ can be obtained from some caterpillar $T'$ by subdividing some
of its leaf edges (not subdividing the spine of $T'$). We find a copy of $T'$ with all its vertices in different blocks, 
and grow each leaf
of $T'$ to an appropriately long path within the block that contained the leaf, by using ``spires'', a variant
of the proof of Gy\'arf\'as~\cite{gyarfas} showing the $\chi$-boundedness of the graphs not containing a fixed path.

Let us turn to the details.
Throughout the remainder of this section,
$\epsilon>0$ is some real number that will be specified later, and
$G$ is a graph with a mass $\mu$, satisfying:
\begin{enumerate}[\indent (1)]
\item $\mu(\{v\})< \epsilon$ for every vertex $v$;
\item $\mu(N(v))< \epsilon$ for every vertex $v$;
and
\item there do not exist two subsets $A,B$ of $V(G)$, anticomplete, with $\mu(A),\mu(B)\ge \epsilon$.
\end{enumerate}
We will show that, for every caterpillar subdivision $T$, if $\epsilon$ is sufficiently small, then $G$ contains $T$ as an induced subgraph, which
will prove \ref{measureversion}. We refer to the three statements above as the ``axioms''.

\begin{thm}\label{components}
Let $X\subseteq V(G)$. If $\mu(X)\ge 3\epsilon$ then $\mu(X')> \mu(X)-\epsilon$ for the vertex set $X'$ of 
some component of $G[X]$.
\end{thm}
\Proof Let the vertex sets of the components of $G[X]$ be $X_1\ll X_k$ say. Choose
$i\ge 1$ minimal such that $\mu(X_1\cup\cdots\cup X_i)\ge \epsilon$. Then from axiom (3), $\mu(X_{i+1}\cup \cdots\cup X_n)<\epsilon$; 
and from the minimality of $i$, $\mu(X_1\cup\cdots\cup X_{i-1})< \epsilon$. But 
$$\mu(X_1\cup\cdots\cup X_{i-1}) + \mu(X_i) + \mu(X_{i+1}\cup \cdots\cup X_n)\ge \mu(X)\ge 3\epsilon,$$
and so $\mu(X_i)\ge \epsilon$. From axiom (3), the union of all other components has mass
less than $\epsilon$, and so $\mu(X_i)> \mu(X)-\epsilon$.
This proves \ref{components}.~\bbox

We observe that since the union of all components of $G[X]$ different from $X'$ has mass less than $\epsilon$,
the set $X'$ in \ref{components} is unique, and we call it the 
{\em big piece} of $X$.

\begin{thm}\label{nested}
Let $X\subseteq Y\subseteq V(G)$. If $\mu(X)\ge 3\epsilon$ then the big piece of $X$ is a subset of the big piece of $Y$.
\end{thm}
\Proof
The big piece of $X$ has mass at least  $\epsilon$, and is a subset of the vertex set of some component of $G[Y]$;
and therefore is a subset of the big piece of $Y$. This proves \ref{nested}.~\bbox

Let $\tau\ge 3$ be an integer.
If $X\subseteq V(G)$, a {\em $\tau$-spire} in $X$ is a sequence $(x_1\ll x_{\tau}, Z)$, where
\begin{itemize}
\item $x_1\ll x_{\tau}$ are the vertices in order of an induced $\tau$-vertex path of $G[X]$;
\item $Z\subseteq X\setminus \{x_1\ll x_{\tau-1}\}$, and $x_{\tau}\in Z$;
\item $x_1\ll x_{\tau-1}$ have no neighbours in $Z\setminus \{x_{\tau}\}$; and
\item $G[Z]$ is connected.
\end{itemize}

\begin{thm}\label{getspire}
Let $\tau\ge 3$ be an integer, and let $X\subseteq V(G)$ with $\mu(X)\ge (\tau+2)\epsilon$;
then there is a $\tau$-spire $(x_1\ll x_{\tau}, Z)$ in $X$ where $\mu(Z)\ge \mu(X)-\tau\epsilon$.
\end{thm}
\Proof
Let $Z_1$ be the big piece of $X$, and choose $x_1\in Z_1$.
Let $Z_2$ be the big piece of $X\setminus N(x_1)$. Since $\mu(X\setminus N(x_1))\ge 3\epsilon$,
from axiom (2) and since $\tau\ge 3$,
\ref{nested} implies that
$Z_2\subseteq Z_1$. Now $x_1$ is a one-vertex component of $G[X\setminus N(x_1)]$, and therefore not its big piece,
by axiom (1);
and since $Z_2\subseteq Z_1$, some neighbour $x_2$ of $x_1$ has a neighbour in $Z_2$.

Inductively, suppose that $2\le i<\tau$, and we have defined $x_1\ll x_i$ and $Z_{i}$, where
\begin{itemize}
\item $x_1\ll x_{i}$ are the vertices in order of an induced $i$-vertex path of $G[X]$;
\item $Z_{i}$ is the big piece of $X\setminus \bigcup_{1\le h\le i-1} N(x_h)$; and
\item $x_i$ has a neighbour in $Z_{i}$.
\end{itemize}
Let $Y_{i+1} = X\setminus \bigcup_{1\le h\le i} N(x_h)$. 
Axiom (2) implies that $\mu(Y_{i+1})\ge \mu(X)-i\epsilon\ge 3\epsilon$.
Let $Z_{i+1}$ be the big piece of $Y_{i+1}$. By \ref{nested}, $Z_{i+1}\subseteq Z_i$, and so some neighbour $x_{i+1}$
of $x_i$ has a neighbour in $Z_{i+1}$.
This completes the inductive definition.

Then $(x_1\ll x_{\tau}, Z_{\tau}\cup \{x_{\tau}\})$ is a $\tau$-spire in $X$, and $\mu(Z_{\tau})\ge \mu(X)-\tau\epsilon$, by 
axiom (2) and \ref{components}.
This proves \ref{getspire}.~\bbox

Let $H$ be a caterpillar, and choose a vertex $v$ which is an end of some path $P$ of $H$ that contains all 
vertices with degree at least two; and
call $v$ the {\em head} of the caterpillar. 
The {\em spine} is the minimal path of $H$ with one end $v$ that contains all vertices of degree at least two.
The pair $(H,v)$ is thus a rooted tree rather than a tree, but we will normally speak of it as a tree 
and let the head be implicit.

Again, let $\tau\ge 3$ be an integer.
A caterpillar is a {\em $\tau$-chrysalis} if
\begin{itemize}
\item its spine has at most $\tau+1$ vertices;
\item every vertex of the spine different from the head has degree exactly $\tau$; and
\item the head has degree at most $\tau-1$, and the head has degree one if the spine has $\tau+1$ vertices.
\end{itemize}
The $\tau$-chrysalis with most vertices therefore has $\tau^2-\tau+2$ vertices, and is unique; let us call it the {\em $\tau$-butterfly}.
It is the only $\tau$-chrysalis in which the spine has $\tau+1$ vertices.

Now let $N$ be a disjoint union of $\tau$-chrysalises $H_1\ll H_k$; we call $N$ a {\em $\tau$-nursery}. 
We define 
$$\phi(N)= \sum_{1\le i\le k}2^{|V(H_i)|}.$$
If $N, M$ are $\tau$-nurseries, 
we say that $M$ is an {\em improvement} of $N$ if 
$M$ has fewer components than $N$ and 
$\phi(M)\ge \phi(N)$.

Returning to the graph $G$ with mass $\mu$, we need to define what it is for a $\tau$-nursery $N$ to be 
``realizable'' in $G$. Let us direct all the edges of $N$ towards the heads; thus, 
for every edge $uv$ of $N$, if $v$ is on the path between $u$ and the head of the component of $N$ containing $u$,
we direct the edge $uv$ from $u$ to $v$. A vertex $v$ of $N$ is a {\em leaf} if it has indegree zero and outdegree one in $N$;
that is, if and only if it does not belong to the spine of its component.
Let $0\le \kappa\le 1$, and for each vertex $v\in V(N)$, let $X_v\subseteq V(G)$, satisfying the following conditions:
\begin{itemize}
\item the sets $X_v\;(v\in V(N))$ are pairwise disjoint;
\item for each leaf $v$ of $N$
there is a $\tau$-spire $(x^1_v\ll x^{\tau}_v, Z_v)$ in $X_v$, and $X_v=\{x^1_v\ll x^{\tau}_v\}\cup Z_v$;
\item for all distinct $u,v\in V(N)$, if $v$ is a leaf then $\{x^1_v\ll x^{\tau}_v\}$ is anticomplete to $X_u$;
\item for all distinct $u,v\in V(N)$, if there is an edge of $G$ between $X_u, X_v$ then either $u,v$ are adjacent in $N$ or both
$u,v$ are heads of components of $N$;
\item for every directed edge $u\rightarrow v$ of $N$, 
$X_u$ covers $X_v$;
\item for each $v\in V(N)$, if $v$ is the head of a component of $N$ then $\mu(X_v)\ge \kappa$.
\end{itemize}
If such a function $X_v\;(v\in V(N))$ exists we call it a {\em $\kappa$-realization} of $N$ in $G$, and say $N$ is 
{\em $\kappa$-realizable} in $G$.
We need:

\begin{thm}\label{stepup}
Let $\tau\ge 3$ be an integer, and let $0\le \kappa,\kappa'\le 1$, with $\kappa\ge 2\kappa' + (\tau+2)\epsilon$. 
Let $N$ be a $\tau$-nursery with at least two components, and in which no component is the $\tau$-butterfly.
If $N$ is $\kappa$-realizable in $G$, there is an improvement $N'$ of $N$ that is
$\kappa'$-realizable in $G$.
\end{thm}
\Proof Let the components of $N$ be $H_1\ll H_k$, where $|V(H_1)|\le \cdots \le |V(H_k)|$, 
and for $1\le i\le k$ let $h_i$ be the head of $H_i$.
Let $X_v\;(v\in V(N))$ be a $\kappa$-realization of $N$ in $G$.
If there exists $i\in \{1\ll k\}$ such that $h_i$ has degree $\tau-1$, choose such a value of $i$, maximum; 
and otherwise, let $i = 1$.
By \ref{getspire}, since $\mu(X_{h_i})\ge \kappa\ge (\tau+2)\epsilon$,
there is a $\tau$-spire $(x_1\ll x_{\tau}, Z)$ say in $X_{h_i}$ where $\mu(Z)\ge \mu(X_{h_i})-\tau\epsilon\ge \epsilon$.

For each $j\in \{1\ll k\}$ with $j\ne i$, let $Y_{h_j}\subseteq X_{h_j}$ be the set of vertices in $X_{h_j}$
with no neighbour in $\{x_1\ll x_{\tau}\}$. Thus 
$$\mu(Y_{h_j})\ge \mu(X_{h_j})-\tau\epsilon\ge \kappa-\tau\epsilon\ge 2(\kappa'+\epsilon)$$
from axiom (2).
Since $G[Z]$ is connected and $x_{\tau}\in Z$, we can number the vertices of $Z$ as $z_1\ll z_n$ say, such that
$z_1 = x_{\tau}$ and $G[\{z_1\ll z_m\}]$ is connected for $1\le m\le n$. 
Since $k\ge 2$, there exists $j\ne i$ with $1\le j\le k$;
but $\mu(Y_{h_j})\ge 2(\kappa'+\epsilon)$, and 
by axiom (3), the set of vertices in $Y_{h_j}$ with no neighbour in $Z$ has mass less than $\epsilon$.
Consequently we may choose $m$ with $1\le m\le n$, minimum such that for some $j\in \{1\ll k\}\setminus \{i\}$,
the set of vertices in $Y_{h_j}$ with no neighbour in $\{z_1\ll z_m\}$ has mass less than $\kappa'+\epsilon$. 
Since no vertex in $Y_{h_j}$ is adjacent to $z_1$, it follows that $m\ge 2$.

\begin{itemize}
\item If $j<i$, it follows that the degree of $h_i$ in $N$ is exactly $\tau-1$. Let $N'$ be the graph obtained from $N$
by adding the edge $h_ih_j$, and deleting all vertices in $V(H_j)\setminus \{h_j\}$. 
Let $H_i'$ be the component
of $N'$ that contains the edge $h_ih_j$, and let us assign its head to be $h_j$.
Thus $H_i'$ is a $\tau$-chrysalis, and so $N'$ is a $\tau$-nursery. Since $N'$ has $k-1$ components and
$|V(H_i)|\ge |V(H_j)|$ (because $i>j$)
it follows that $\phi(N')\ge \phi(N)$,
and $N'$ is an improvement of $N$.
\item If $j>i$, it follows that the degree of $h_j$ in $N$ is at most $\tau-2$. Let $N'$ be the graph obtained from $N$
by adding the edge $h_ih_j$, and deleting all vertices in $V(H_i)\setminus \{h_i\}$. 
Let $H_j'$ be the component
of $N'$ that contains the edge $h_ih_j$, and let us assign its head to be $h_j$.
Thus $H_j'$ is a $\tau$-chrysalis, 
and again $N'$ is an improvement of $N$.
\end{itemize}
For each $v\in V(N')$ define $X_v'$ as follows: 
\begin{itemize}
\item if $v\ne \{h_1\ll h_k\}$ let $X_v'=X_v$;
\item let $X_{h_i}'=\{z_1\ll z_m\}\cup \{x_1\ll x_{\tau}\}$;
\item let $X_{h_j}'$ be the set of vertices in $Y_{h_j}$ with a neighbour in $\{z_1\ll z_m\}$;
\item for $1\le \ell\le k$ with $\ell\ne i,j$, let $X_{h_{\ell}}'$ be the set of vertices in $Y_{h_{\ell}}$ with no neighbour in
$\{z_1\ll z_m\}$.
\end{itemize}
We see that $X_{h_i}'$ covers $X_{h_j}'$, and has no edges to $X_{h_{\ell}}'$ for $1\le \ell\le k$ with $\ell\ne i,j$.
Moreover, $\mu(X_{h_j}')\ge \kappa'$. Let $1\le \ell\le k$ with $\ell\ne i,j$; then, since $m\ge 2$ and
from the choice of $m$,
the mass of the set of vertices in $Y_{h_{\ell}}$ with no neighbour in $\{z_1\ll z_{m-1}\}$ is at least $\kappa'+\epsilon$.
Hence $\mu(X_{h_{\ell}}')\ge \kappa'$. It follows that the function $X_v'\;(v\in V(N'))$
is a $\kappa'$-realization of $N'$ in $G$. This proves \ref{stepup}.~\bbox

Now let $T$ be a caterpillar subdivision. We say that an integer $\tau\ge 3$ {\em fits $T$} if
\begin{itemize} 
\item there is a path of $T$ with at most $\tau$ vertices containing all vertices of $T$ of degree more than two; 
\item $T$ has maximum degree at most $\tau$; and 
\item every path of $T$ in which every internal vertex has degree two in $T$ has at most $\tau$ vertices.  
\end{itemize}

\begin{thm}\label{endstop}
Let $T$ be a caterpillar subdivision, and let $\tau$ fit $T$. If $G$ is $T$-free then 
for $\kappa>0$, the $\tau$-butterfly is not $\kappa$-realizable in $G$.
\end{thm}
\Proof
Suppose that $X_v\;(v\in V(N))$ is a $\kappa$-realization in $G$ of the $\tau$-butterfly $N$. Now $N$ is connected, and 
since $|V(N)|=\tau^2-\tau+2$, the spine
of $N$ has exactly $\tau+1$ vertices and they all have degree $\tau$ except the head which has degree one. 
Let the spine of $N$ have vertices $v_0,v_1\ll v_{\tau}$
in order, where $v_0$ is the head of $N$. Since $\mu(X_{v_0})\ge \kappa>0$, it follows that $X_{v_0}\ne \emptyset$;
choose $p_{v_0}\in X_{v_0}$. For $1\le i\le \tau$, choose $p_{v_i}\in X_{v_i}$ adjacent to $p_{v_{i-1}}$; this is 
possible since $X_{v_i}$ covers $X_{v_{i-1}}$. Now let $u$ be a leaf of $N$, with neighbour $v$ say.
From the definition of a realization,
there is a $\tau$-spire $(x^1_u\ll x^{\tau}_u, Z_u)$ in $X_u$, and $X_u=\{x^1_u\ll x^{\tau}_u\}\cup Z_u$. Since $p_v$
has a neighbour in $X_u$, and $G[Z_u]$ is connected and contains $x^{\tau}_u$, and none of $x^1_u\ll x^{\tau-1}_u$
have neighbours in $Z_u\setminus \{x^{\tau}_u\}$, there is an induced path $P_u$ 
with $\tau$ vertices, with one end $p_v$ and with all other vertices in $X_u$. Let $H$ be the induced
subgraph of $G$ consisting of
the union of all these paths $P_u$ (over all leaves $u$ of $N$) and the path induced on $\{p_{v_0}\ll p_{v_{\tau}}\}$;
then $T$ is isomorphic to an induced subgraph of $H$, contradicting that $G$ is $T$-free. This proves \ref{endstop}.~\bbox

Now we can prove the main theorem \ref{measureversion}, which we restate.
\begin{thm}\label{measureversion2}
For every caterpillar subdivision $T$, there exists $\epsilon>0$ such that for every $T$-free graph $G$, and
mass $\mu$ on $G$, either
\begin{itemize}
\item $\mu(\{v\})\ge \epsilon$ for some vertex $v$; or
\item $\mu(N(v))\ge \epsilon$ for some vertex $v$; or
\item there are two anticomplete sets of vertices $A,B$, where $\mu(A), \mu(B)\ge \epsilon$.
\end{itemize}
\end{thm}
\Proof
Choose $\tau$ fitting $T$,
and let $p = 2^{\tau^2}$. Define $\epsilon$ such that $\epsilon^{-1} = p2^{p}(\tau+3)$. 
We will show that $\epsilon$ satisfies the theorem.
Suppose not, and choose a $T$-free graph $G$, and mass $\mu$ on $G$ not satisfying the theorem (and therefore satisfying
the axioms).
For $0\le i\le p$ define 
$\kappa_i = 2^{-i}p^{-1} - (\tau+2)\epsilon$.
Thus $0\le \kappa_i\le 1$ for each $i$. Moreover,
$\kappa_{p}=\epsilon$, and for $1\le i\le p$,
$$\kappa_{i-1}= 2\kappa_i  + (\tau+2)\epsilon.$$
Choose $X_1\ll X_P\subseteq V(G)$, pairwise disjoint, with $\kappa_0 \le \mu(X_i)< \kappa_0+\epsilon$ for $1\le i\le P$, 
with $P$ maximum.
We claim that $P\ge p$; for suppose not. Then the union of $X_1\ll X_P$ has mass at most $(p-1)(\kappa_0+\epsilon)$, 
and since
$(p-1)(\kappa_0+\epsilon)\le 1- \kappa_0$, there exists a set of mass at least $\kappa_0$ disjoint from this union. Choose
such a set, $X_{P+1}$ say, minimal; then from the minimality of $X_{P+1}$, and 
since $\mu(\{v\})< \epsilon$ for each vertex $v$, it follows that $\mu(X_{P+1})< \kappa_0+\epsilon$,
contrary to the maximality of $P$. This proves that $P\ge p$.

Let $N_0$ be the $\tau$-nursery with $p$ components, each an isolated vertex.
It follows that $N_0$ is $\kappa_0$-realizable in $G$ and $\phi(N_0) = 2p$.
Choose a sequence $N_1\ll N_q$ of $\tau$-nurseries, such that for $1\le i\le q$, $N_i$ is an improvement of $N_{i-1}$, and
$N_i$ is $\kappa_i$-realizable in $G$, with $q$ maximum. It follows that 
$\phi(N_i)\ge \phi(N_{i-1})$ for $1\le i\le q$, from the definition of an improvement, and so $\phi(N_q)\ge 2p$, and
in particular, $N_q$ is nonnull. But $N_i$ has at most $p-i$ components for $0\le i\le q$,
and so $q\le p-1$. Thus $\kappa_{q+1}$ is defined.
By \ref{endstop} no component of $N_q$ is the $\tau$-butterfly, and so $N_q$ has at most one component by
\ref{stepup}, and therefore has at most $\tau^2-\tau+1$ vertices. 
But $\phi(N_q) \ge 2p$, which is impossible.

Thus there is no such pair $G,\mu$. This proves \ref{measureversion2}.~\bbox


\begin{thebibliography}{99}
\bibitem{lagoutte} N. Bousquet, A. Lagoutte, and S. Thomass\'e, ``The Erd\H{o}s-Hajnal conjecture for paths and
antipaths'', {\em J. Combinatorial Theory, Ser. B}, {\bf 113} (2015), 261--264.
\bibitem{hooks} K. Choromanski, D. Falik, A. Liebenau, V. Patel, and M. Pilipczuk, ``Excluding hooks and their complements'', {\em Electron.~J.~Combin.}, {\bf 25(3)} (2018), \#P3.27.
\bibitem{maria} M. Chudnovsky, ``The Erd\H{o}s-Hajnal conjecture - a survey'', {\em J. Graph Theory}
{\bf 75} (2014), 178--190.
\bibitem{trees} M. Chudnovsky, A. Scott, P. Seymour and S. Spirkl, ``Trees and linear anticomplete pairs'',
submitted for publication.
\bibitem{girth} P. Erd\H{o}s, ``Graph theory and probability'', {\em Canadian J. Math.} {\bf 11} (1959), 34--38.
\bibitem{EH0} P. Erd\H{o}s and A. Hajnal, ``On spanned subgraphs of graphs'', 
{\em Graphentheorie und Ihre Anwendungen} (Oberhof, 1977), \verb++{www.renyi.hu/\raisebox{-1ex}{\textasciitilde}p\_erdos/1977-19.pdf}.
\bibitem{EH}  P. Erd\H{o}s and A. Hajnal, ``Ramsey-type theorems'',
{\em  Discrete Applied Mathematics} {\bf 25} (1989), 37--52.
\bibitem{fox0} J. Fox and B. Sudakov, ``Induced Ramsey-type theorems'',
{\em Advances in Math.} {\bf 219} (2008), 1771--1800.
\bibitem{fox} J. Fox, B. Sudakov and A. Suk, ``Erd\H{o}s-Hajnal conjecture for graphs with bounded VC-dimension'',
{\tt arXiv:1710.03745}.
\bibitem{gyarfastree}
A. Gy\'arf\'as, ``On Ramsey covering-numbers'',
{\em Coll. Math. Soc. J\'anos Bolyai}, in {\em Infinite and Finite Sets},
North Holland/American Elsevier, New York (1975), {\bf 10}.
\bibitem{gyarfas}
A. Gy\'arf\'as, ``Problems from the world surrounding perfect graphs'', {\em Proceedings of
the International Conference on Combinatorial Analysis and its Applications},  (Pokrzywna, 1985),
{\em Zastos. Mat.}  {\bf 19} (1987), 413--441.
\bibitem{liebenau} A. Liebenau and M. Pilipczuk, ``The Erd\H{o}s-Hajnal conjecture for caterpillars and their complements'',
{\tt arXiv:1710.08701}.
\bibitem{rodl} V. R\"odl, ``On universality of graphs with uniformly distributed edges'', 
{\em Discrete Math.} {\bf 59} (1986), 125--134.
\bibitem{sumner}
D.P. Sumner, ``Subtrees of a graph and chromatic number'', in
{\em The Theory and Applications of Graphs}, (G. Chartrand, ed.),
John Wiley \& Sons, New York (1981), 557--576.
\end{thebibliography}
\end{document}